\documentclass[oneside,notitlepage,12pt]{article}

\usepackage{amssymb}
\usepackage[leqno]{amsmath}
\usepackage{amsfonts}
\usepackage{amsopn}
\usepackage{amstext}
\usepackage{amsthm}

\usepackage{tikz}
\usetikzlibrary{cd}

\usepackage{enumitem}

\usepackage{verbatim}
\usepackage[colorlinks, backref]{hyperref}
\usepackage{makeidx}

\usepackage{calrsfs}
\usepackage{clock} 

\providecommand{\cal}{\mathcal}

\newenvironment{pf}{\begin{proof}}{\end{proof}}

\newcommand{\Aaa}{{\cal{A}}}

\newcommand{\Cee}{{\cal{C}}}

\newcommand{\Ef}{{\cal{F}}}

\newcommand{\Z}{{\mathbb{Z}}} 
\newcommand{\Q}{{\mathbb{Q}}} 

\newcommand{\al}{\alpha}

\newcommand{\sig}{\sigma}

\renewcommand{\phi}{\varphi}
\renewcommand{\rho}{\varrho}

\newcommand{\rest}{\restriction}

\newcommand{\ntr}{{n\in\omega}}

\newcommand{\loe}{\leq}
\newcommand{\goe}{\geq}

\newcommand{\subs}{\subseteq}

\renewcommand{\iff}{\Longleftrightarrow}

\newcommand{\id}[1]{{\operatorname{i\!d}_{#1}}} 

\newcommand{\by}[1]{/{#1}}

\newtheorem{tw}{Theorem}[section]
\newtheorem{wn}[tw]{Corollary}
\newtheorem{lm}[tw]{Lemma}
\newtheorem{prop}[tw]{Proposition}
\newtheorem{claim}[tw]{Claim}
\theoremstyle{definition}

\newtheorem{question}[tw]{Question}
\theoremstyle{remark}

\newcommand{\set}[1]{\{#1\}}
\newcommand{\setof}[2]{\{#1\colon #2\}}

\newcommand{\seq}[1]{\langle #1 \rangle}

\newcommand{\dn}[2]{\{#1,#2\}} 
\newcommand{\pair}[2]{\langle #1, #2 \rangle} 
\newcommand{\map}[3]{#1\colon #2 \to #3} 

\newcommand{\fra}{Fra\"iss\'e}

\providecommand{\nat}{\omega}

\newcommand{\aut}{\operatorname{Aut}}

\newcommand{\iso}{\approx}

\newcommand{\emb}[1]{{{}^{\operatorname{emb}}{#1}}}

\newcommand{\abs}[1]{|#1|}

\newcommand{\cmp}{\circ} 

\providecommand{\ar}{\arrow}

\newcommand{\Katetov}{Kat\v{e}tov}
\newcommand{\Ens}{\mathfrak{S}}

\newcommand{\age}{\operatorname{age}}

\title{Uniform homogeneity}
\author{
{\sc Wies{\l}aw Kubi\'s}\\[1mm]
{\small Institute of Mathematics, Czech Academy of Sciences (CZECHIA)}\\[2mm]
\and
{\sc Bori\v{s}a Kuzeljevi\'c}\\[1mm]
{\small Institute of Mathematics, Czech Academy of Sciences (CZECHIA)}\\
{\small Department of Mathematics and Informatics, University of Novi Sad (SERBIA)}
}

\date{\clocktime\today}

\makeindex

\begin{document}

\maketitle

\begin{abstract}
We discuss some finite homogeneous structures, addressing the question of universality of their automorphism groups.
We also study the existence of so-called \Katetov\ functors in finite categories of embeddings or homomorphisms.

\end{abstract}

\tableofcontents

\section{Introduction}

By a \emph{structure} we mean a model of a fixed first-order language, possibly involving algebraic operations. The notions of a substructure, embedding, homomorphisms, etc. are defined in the obvious way.
A structure $M$ is \emph{homogeneous} if every isomorphism between finitely generated substructures of $M$ extends to an automorphism of $M$. Some authors call it \emph{ultra-homogeneity} in order to distinguish it from point-homogeneity, where \emph{points} are one-element substructures. In the presence of algebraic operations or unary predicates, either points may generate infinite substructures or there may be several types of points, therefore point-homogeneity does not say that the automorphism group acts transitively.

We are interested in finite homogeneous structures. Perhaps the simplest ones are just sets, namely, where the language is empty. Automorphisms are bijections, while embeddings are one-to-one mappings.
Purely algebraic examples are finite cyclic groups.
Note that the infinite cyclic group $\Z$ is not homogeneous, as there is no automorphism mapping $k\Z$ onto $\ell\Z$ whenever $k,\ell>0$ are distinct.

We are interested in the following natural question: Given a (possibly finite) homogeneous structure $M$, is its automorphism group $\aut(M)$ \emph{universal} in the sense that it embeds all the groups of the form $\aut(X)$ with $X$ a substructure of $M$?
Positive results are usually achieved with the help of a functor from the category of all isomorphisms between substructures of $M$ into the group of automorphisms of $M$ (recall that each monoid $G$ is a category with the single object $G$, whose arrows are the elements of $G$).

Let us call a structure $M$ \emph{uniformly homogeneous} if there is a functor $K$ from the category of all isomorphisms between finitely generated substructures of $M$ into the group $\aut(M)$ such that $K(f)$ is an extension of $f$ for each $f$.
In other words, $K$ is an extension operator on isomorphisms of finitely generated substructures of $M$ satisfying
\begin{enumerate}[itemsep=0pt]
	\item[(1)] $K(\id{A}) = \id M$,
	\item[(2)] $K(f) \in \aut(M)$ and $K(f)$ extends $f$,
	\item[(3)] $K(f \cmp g) = K(f) \cmp K(g)$,
\end{enumerate}
for every isomorphisms $\map f A B$, $\map g C A$ with $A,B,C$ finitely generated substructures of $M$.
Clearly, every uniformly homogeneous structure is homogeneous. Most of the well known infinite homogeneous structures are uniformly homogeneous, as it is argued in~\cite{KubMas}.
Proposition~\ref{PropSDnldvno} below gives a useful criterion for uniform homogeneity.

We are particularly interested in countable homogeneous structures whose age (that is, the class of all finitely generated substructures) is finite. Of course, finite homogeneous structures have this property, however there exist also infinite ones.

\subsection{Notation}

We shall use standard notation concerning model theory and set theory.
In particular, $n = \{0,1,\dots, n-1\}$ for every positive integer $n$. The set of natural numbers is $\nat = \{0,1,\cdots\}$.
If $f:X\to Y$ is a function from $X$ to $Y$, and $Z\subseteq X$, then $f[Z]$ denotes the image of $Z$ under $f$, i.e. $f[Z]=\set{f(x):x\in Z}$.
The greatest common divisor of two numbers $n,m$ will be denoted, as usual, by $\gcd(n,m)$.
As we said before, a \emph{structure} is a model of a countable first-order language.
The \emph{age} of a structure $M$ will be denoted by $\age(M)$. Recall that $\age(M)$ is the class of all finitely generated models embeddable into $M$.
By an \emph{embedding} we mean a one-to-one mapping that is an isomorphism onto its image. A \emph{substructure} of $M$ is a subset closed under all operations of $M$ (in particular containing all constants) with induced relations. Namely, given an $n$-ary relation $R$ in $M$, the induced relation in $X \subs M$ is $R^X$ defined by $R^X(a_0,\dots, a_{n-1})$ if and only if $M \models R(a_0, \dots, a_{n-1})$ for every $a_0, \dots, a_{n-1} \in X$.
We shall write $X \loe M$ to say that $X$ is a substructure of $M$.
Recall that a structure $M$ is \emph{homogeneous} if every isomorphism $\map f X Y$ between finitely generated substructures of $M$ extends to an automorphism of $M$.
In case $M$ is countable (as we will always assume), this is equivalent to the \emph{extension property} of $M$ saying that for every embeddings $\map e A M$, $\map f A B$, with $A,B \in \age(M)$, there exists an embedding $\map g B M$ satisfying $e = g \cmp f$.
For details on \fra\ theory we refer to Chapter 7 of Hodges' monograph~\cite{Hodges} (see also the original paper of \fra~\cite{Fraisse}).
Another relevant notion is of set-homogeneity.
A structure $M$ is \emph{set-homogeneous} if for every two isomorphic finitely generated substructures $A$ and $B$ of $M$, there is an automorphism $f$ of $M$ such that $f[A]=B$.
Note that every homogeneous structure is also set-homogeneous.

\subsection{A characterization of uniform homogeneity}

\begin{prop}\label{PropSDnldvno}
	Let $M$ be a set-homogeneous structure. Then $M$ is uniformly homogeneous if and only if for every finitely generated substructure $A$ of $M$ there exists a homomorphism $\map{E_A}{\aut(A)}{\aut(M)}$ such that $E_A(h)$ is an extension of $h$ for every $h \in \aut(A)$.
\end{prop}

\begin{pf}
	Clearly, the condition is necessary, as we may set $E_A = K \rest \aut(A)$.
	
	In order to show sufficiency, note that it is enough to define $K$ satisfying (1), (2) and (3) for each isomorphism class separately.
	Fix a finitely generated $A \loe M$.
	Let $\Aaa$ be the family of all substructures of $M$ isomorphic to $A$.
	For each $X \in \Aaa$ choose $\phi_X \in \aut(M)$ such that $\phi_X \rest A$ is an isomorphism onto $X$. Here we used the set-homogeneity of $M$.
	Now, given an isomorphism $\map f X Y$ with $X, Y \in \Aaa$, define
	$$K(f) = \phi_Y \cmp E_A(\phi^{-1}_Y\cmp f \cmp \phi_X \rest A) \cmp \phi^{-1}_X.$$
	Note that $\phi^{-1}_Y\cmp f \cmp \phi_X \rest A \in \aut(A)$, therefore $K$ is well defined.
	Given $x \in X$, we have
	\begin{align*}
	K(f)(x) &= \phi_Y \cmp E_A(\phi^{-1}_Y\cmp f \cmp \phi_X \rest A) \cmp \phi^{-1}_X(x) \\
	&= \phi_Y \cmp \phi^{-1}_Y\cmp f \cmp \phi_X \cmp \phi^{-1}_X(x) = f(x),	
	\end{align*}
	therefore $K(f)$ extends $f$. Clearly, $K(f) \in \aut(M)$. It is also clear that $K(\id X) = \id M$.
	Finally, given an isomorphism $\map g Y Z$, we have
	\begin{align*}
		K(g) \cmp K(f) &= (\phi_Z \cmp E_A(\phi^{-1}_Z\cmp g \cmp \phi_Y \rest A) \cmp \phi^{-1}_Y) \\ &\cmp (\phi_Y \cmp E_A(\phi^{-1}_Y\cmp f \cmp \phi_X \rest A) \cmp \phi^{-1}_X) \\
		&= \phi_Z \cmp E_A(\phi^{-1}_Z \cmp g \cmp \phi_Y \rest A \cmp \phi^{-1}_Y \cmp f \cmp \phi_X \rest A) \cmp \phi^{-1}_X \\
		&= \phi_Z \cmp E_A(\phi^{-1}_Z \cmp g \cmp f \cmp \phi_X \rest A) \cmp \phi^{-1}_X = K(g \cmp f).
	\end{align*}
	This completes the proof.
\end{pf}

\subsection{\Katetov\ functors}

Let $\Ef$ be a \fra\ class in a countable signature. We denote by $\emb{\Ef}$ the category of all embeddings between structures from $\Ef$.
Let $\sig \Ef$ denote the class of all countable structures whose age is contained in $\Ef$.
Following~\cite{KubMas}, a \emph{\Katetov\ functor} on $\Ef$ is a functor $\map K{\emb{\Ef}}{\emb{\sig\Ef}}$ for which there exists a natural transformation $\eta$ from the identity to $K$, such that for every $X \in \Ef$, for every one-point extension $\map e X {X'}$ there exists an embedding $\map f {X'}{K(X)}$ satisfying $\eta = f \cmp e$.
Note that this is automatically satisfied whenever $K(X)$ is the \fra\ limit of $\Ef$.
A \Katetov\ functor $K$ will be called \emph{ultimate} if $K(X)$ is (isomorphic to) the \fra\ limit of $\Ef$ for every $X \in \Ef$.
One of the basic results in~\cite{KubMas} says that every \Katetov\ functor extends to $\sig \Ef$ and its $\omega$th iteration is an ultimate \Katetov\ functor.
Thus, we may assume that every \Katetov\ functor takes values in the monoid $\emb{U}$ of all embeddings from $U$ to $U$, where $U$ is the \fra\ limit of $\Ef$.
Clearly, if $\Ef$ admits a \Katetov\ functor then its \fra\ limit is uniformly homogeneous.
The converse totally fails, at least for finite \fra\ classes.

\begin{prop}
	Let $\Ef$ be a \fra\ class whose \fra\ limit $U$ is finite.
	Assume $A$ is a substructure of $U$ such that there exists $h \in \aut(U)$ with $h \ne \id U$ and $h \rest A = \id A$.
	Then $\Ef$ admits no \Katetov\ functor.
\end{prop}

\begin{pf}
	Suppose $K$ is a \Katetov\ functor on $\Ef$. By the remarks above, we may assume $K(A) = U = K(U)$.
	Let $\eta$ be the associated natural transformation.
	Let $\map i A U$ be the inclusion.
	Then $h \cmp i = i$, therefore $K(h) \cmp K(i) = K(h \cmp i) = K(i)$.
	On the other hand, $K(i)$ is an embedding of $U$ into itself, therefore it is an isomorphism, because $U$ is finite.
	Hence $K(h) = \id U$.
	Now, using the fact that $\eta$ is a natural transformation from the identity functor to $K$, we obtain
	$$\eta_U \cmp h = K(h) \cmp \eta_U = \id U \cmp \eta_U = \eta_U,$$
	therefore $h = \id U$, because $\eta_U$ is an embedding.
\end{pf}

\begin{wn}
	Let $\Ef$ be the class of all sets of cardinality $\loe n$, where $n \goe 3$. Then $\Ef$ is a \fra\ class with no \Katetov\ functor.
\end{wn}

Given $\ntr$, let $\Ens(n)$ denote the class of all sets of cardinality $\loe n$.
Without loss of generality, we may assume that $\Ens(n)$ consists of subsets of $n = \{0,1, \dots, n-1\}$.

Given a bijection $\map f A B$, define $K(A)=K(B)=n$ and $\map{K(f)}{n}{n}$ in such a way that the set $n \setminus A$ is mapped in a strictly increasing way onto the set $n \setminus B$.
It is rather clear that $K(g \cmp f) = K(g) \cmp K(f)$, because we deal with bijections. Also, $K(\id A) = \id n$. Thus, $K$ is a functor. This shows that every finite set is uniformly homogeneous.

It is well known and very easy to verify that the group of permutations $S_n = \aut(n)$ is universal for the class $\setof{S_k}{k \loe n}$. The embedding of $S_k$ into $S_n$ is given by $h \mapsto h \cup \id{n \setminus k}$.

\section{A simple digraph with six vertices}

A \emph{simple digraph} is a structure of the form $\pair X E$, where $E$ is a binary relation on $X$. The elements of $X$ are usually called \emph{vertices}, while the elements of $E$ are called \emph{arrows} (some authors call them \emph{edges}).

Our goal is to describe a homogeneous simple digraph with $6$ vertices, with no \Katetov\ functor on the category of isomorphisms between its substructures.
Our graph is actually described in the following picture.
$$
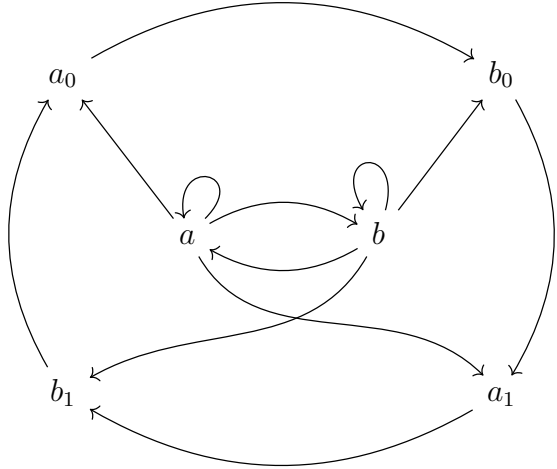
\begin{tikzcd}
a_0 \ar[rrrr, bend left] & & & & b_0 \ar[dddd, bend left] \\
& & & & \\
& a \ar[luu] \ar[rr, bend left] \ar[in=100,out=45,loop] \ar[out=300,ddrrr] & & b \ar[ruu] \ar[in=30,,out=240,llldd] \ar[ll, bend left] \ar[in=120,out=75,loop] & \\
& & & & \\
b_1 \ar[uuuu, bend left] & & & & a_1 \ar[llll, bend left]
\end{tikzcd}$$
Let us denote this digraph by $M$.
Formally, the universe of $M$ is $\{a,b,a_0,b_0,a_1,b_1\}$ and the relation is
\begin{align*}
	E = &\{\pair aa, \pair ab, \pair ba, \pair bb, \pair{a_0}{b_0},\pair{b_0}{a_1},\pair{a_1}{b_1},\pair{b_1}{a_0},\\
	&\pair{a}{a_0},\pair{a}{a_1},\pair{b}{b_0},\pair{b}{b_1}\}.
\end{align*}
Note that $a,b$ are the only vertices with loops. Hence, every automorphism of $M$ preserves each of the cycles $\{a,b\}$ and $\{a_0, b_0, a_1, b_1\}$.
Furthermore, every automorphism of $M$ ``rotates" the cycle $\{a_0,b_0,a_1,b_1\}$.
Namely, let $\eta \in \aut(M)$ be such that $\eta(a_0) = b_0$.
Then $\eta(a) = b$, $\eta(b_0) = a_1$, $\eta(a_1) = b_1$, and $\eta(b_1) = a_0$.
The same argument shows that every automorphism of $M$ is determined by its value on $a_0$.
It follows that $\eta$ generates $\aut(M)$ and consequently $\aut(M) \iso \Z_4$.

\begin{tw}
	$M$ is homogeneous but not uniformly homogeneous.
\end{tw}

\begin{pf}
	We first show that $M$ is not uniformly homogeneous.
	Let $h_0 \in \aut(\dn a b)$ be the non-trivial involution and let $h \in \aut(M)$ be its extension.
	If $h(a_0) = b_0$ then $h^2(a_0) = h(b_0) = a_1$, therefore $h^2 \ne \id M$.
	If $h(a_0) = b_1$ then $h^2(a_0) = h(b_1) = a_1$, thus again $h^2 \ne \id M$.
	Hence there is no involution of $M$ extending $h_0$.
	
	It remains to check that $M$ is homogeneous.
	Let $\eta$ be the automorphism introduced in the previous paragraph generating $\aut(M)\iso \Z_4$.
	Thus, $\eta(a_0)=b_0$.
	Denote also $C=\set{a_0,b_0,a_1,b_1}$.
	First we prove a short claim.
	
	\begin{claim}\label{help}
	 Suppose that $\psi\in \aut(M\rest C)$.
	 Then there is $i<4$ so that $\eta^i\rest C=\psi$.
	\end{claim}
	
	\begin{pf}
	 It is clear that $M\rest C$ is isomorphic to the 4-cycle $\overrightarrow{C_4}$, so $\aut(M\rest C)\iso \Z_4$, and $\psi$ is completely determined by $\psi(a_0)$.
	 If $\psi(a_0)=a_0$, then $\psi=\id C$, so $i=0$ satisfies the conclusion of the claim.
	 If $\psi(a_0)=b_0$, then $i=1$ works, and all the other cases are handled similarly.
	\end{pf}

	Take arbitrary non-empty substructures $A$ and $B$ of $M$, and let $\phi:A\to B$ be an isomorphism.
	Notice that for every choice of $A$, $B$, and $\phi$, it must be that $\phi[A\cap C]=B\cap C$ and $\phi[A\setminus C]=B\setminus C$.
	We distinguish two cases, depending on the cardinality of the set $A\setminus\set{a,b}$.
	
	Case 1: $\abs{A\setminus \set{a,b}}=0$.
	In this case there are two possibilities, either $\abs{A}=1$, or $\abs{A}=2$.
	If $\abs{A}=1$, then either $\phi[A]=A$, in which case $\phi$ is the identity, and it can be extended to $\eta^0$, or $\phi[A]\cap A=\emptyset$, so that $\phi$ can be extended to $\eta^1$.
	If $\abs{A}=2$ (i.e. $\phi[A]=A$ and $A=\set{a,b}$), then either $\phi(a)=a$ in which case $\phi(b)=b$ and $\phi$ can be extended to $\eta^0$, or $\phi(a)=b$ in which case $\phi(b)=a$ and $\phi$ can be extended to $\eta^1$.
	
	Case 2: $\abs{A\setminus \set{a,b}}>0$.
	Note that in this case $A\cap C=A\setminus \set{a,b}$.
	Denote $\theta=\phi\rest(A\cap C)$.
	Note that $\theta$ is well defined because $A\cap C$ is non empty.
	Since the 4-cycle $\overrightarrow{C_4}$ is homogeneous \cite{Lachlan}, $M\rest C$ is also homogeneous being isomorphic to it.
	So there is some $\psi\in\aut(M\rest C)$ extending $\theta$.
	By Claim \ref{help}, there is $i<4$ such that $\eta^i\rest C=\psi\rest C$.
	This means that $\eta^i\rest(A\setminus\set{a,b})=\phi\rest(A\setminus\set{a,b})$.
	To finish the proof we have to show that $\eta^i\rest A=\phi$.
	Suppose that this is not the case.
	This means that there is a point $x\in A$ such that $\eta^i(x)\neq \phi(x)$.
	Point $x$ cannot be in $C$ by the choice of $i$, so it must be in $\set{a,b}$.
	Take any point $y\in A\cap C$, and let $z=\phi(y)=\eta^i(y)$.
	Note that $z\in C$, while the assumption $\phi(x)\neq \eta^i(x)$ implies that $\set{\phi(x),\eta^i(x)}=\set{a,b}$.
	So since $\phi$ and $\eta^i$ are isomorphisms, and by the definition of $E$, it must be that
	\[
	 \seq{x,y}\in E \iff \seq{\phi(x),z}\in E \iff \seq{\eta^i(x),z}\notin E \iff \seq{x,y}\notin E,
	\]
	which is clearly impossible.	
	
\end{pf}

\begin{question}
	Does there exist a finite homogeneous structure whose domain has less than six points, and which is not uniformly homogeneous?
\end{question}

\section{Finite cyclic groups}

Most of the results of this section should be well known, however we look at cyclic groups from the perspective of \fra\ theory.
As it happens, every finite cyclic group is uniformly homogeneous.
In fact, homogeneity follows directly from the following easy fact. Uniform homogeneity will follow from the existence of a \Katetov\ functor on the \fra\ class of all finite cyclic groups, whose limit is $\Q \by \Z$.

Every homomorphism of cyclic groups $\map{f}{\Z_m}{\Z_n}$ is determined by $f(1)$, where $1$ is the generator of $\Z_m$. Specifically, $f(i) = f(1) \cdot i$ modulo $n$ for $i < m$. We will write $f = \hat a$, where $a = f(1)$.

\begin{lm}\label{LMwbngfojwneof}
	Let $\map {e,f}{\Z_k}{\Z_n}$ be two embeddings.
	Then there exists an automorphism $\map{h}{\Z_n}{\Z_n}$ such that $f = h \cmp e$.
\end{lm}

\begin{pf}
	Obviously, $n = k \ell$ for some integer $\ell > 0$.
	We may assume that $e = \hat \ell$, which is the canonical embedding.
	Then $f = \hat a$, where $\gcd(a,k \ell) = \ell$. In other words, $a = \ell b$, where $b$ and $k \ell$ are coprime.
	Thus $h = \hat b$ is an automorphism of $\Z_{k\ell}$ and $h \cmp e = f$.
\end{pf}

It follows that $\Z_n$ has the extension property. Its age consists of the trivial group plus all groups of the form $\Z_k$, where $k$ is a divisor of $n$.
Thus, finite cyclic groups are homogeneous. They turn out to be uniformly homogeneous. We shall prove it via a \Katetov\ functor on the class of all finite cyclic groups $\Cee$.
Let us remark that $\Cee$ is hereditary (since subgroups of a cyclic group are cyclic) and has the amalgamation property.
Indeed, given embeddings $\map{f}{\Z_k}{\Z_m}$, $\map{g}{\Z_k}{\Z_n}$, we may replace them by $\map{f' \cmp f}{\Z_k}{\Z_{mn}}$, $\map{g' \cmp g}{\Z_k}{\Z_{mn}}$ and then use Lemma~\ref{LMwbngfojwneof} to get an automorphism $\map{h}{\Z_{mn}}{\Z_{mn}}$ satisfying $\phi \cmp f' \cmp f = g' \cmp g$.
It is not hard to see that the \fra\ limit of $\Cee$ is $\Q \by \Z$ which is isomorphic to the group of all roots of unity in the complex plane.

\begin{tw}
	The class of all finite cyclic groups admits a \Katetov\ functor.
\end{tw}

\begin{pf}
	First let us introduce some notation.
	To shorten the statements we denote $U=\Q \by \Z$, while the set of all prime numbers is denoted $\mathbb P$.
	If $n$ is an integer, and $p$ is a prime, then $[n]_p$ denotes the number $\frac{n}{p^{\al}}$, where $\al$ is a non-negative integer such that $p^{\al}\mid n$ and $p^{\al+1}\nmid n$.
	By Theorem I.8.1 in \cite{Lang}
	\[
	 U\cong \bigoplus_{p\in \mathbb P}U(p),
	\]
	where for a prime $p$, the group $U(p)$ is the subgroup of all elements from $U$ which can be represented by a rational number $\frac{a}{p^{\al}}$ with $a\in \Z$ and $\al$ some positive integer.
	So in this proof, whenever we write $x\in U$ we will assume that $x=\seq{x_p:p\in \mathbb P}$ and that $x_p\in U(p)$ for each $p\in \mathbb P$.
	Notice that if $e:\Z_m\to \Z_{mk}$ is an embedding, then there is a number $n$ such that $n=k\cdot t$ for some t satisfying $\gcd(t,mk)=1$, and that for each $x\in \Z_m$, $e(x)=n\cdot x\ (\operatorname{mod} mk)$.
	Whenever we are given such an embedding $e$, we will assume that we are also given a number $n$ with the mentioned properties, and denote $e$ by $\hat n$.
	For a finite cyclic group $\Z_m$, let $\eta_m:\Z_m\to U$ be defined as follows. 
	For $l<m$ let $\eta_m(l)=\seq{\frac{l\cdot [m]_p}{m}:p\in \mathbb P}$.
	Finally, if $\hat n$ is an embedding between $\Z_m$ and $\Z_{mk}$ define $K(\hat n):U\to U$ so that for $x\in U$, $K(\hat n)(x)=\seq{[n]_p\cdot x_p:p\in \mathbb P}$.
	This is well defined because $\gcd([n]_p,p)=1$.
	To prove that $K$ is a Kat\v etov functor, it is enough to prove that:
	\begin{enumerate}
	 \item if $\hat n: \Z_m\to \Z_{mk}$ is an embedding, then $K(\hat n):U\to U$ is also an embedding;
	 \item if $\hat n_1:\Z_m\to Z_{mk}$ and $\hat n_2:\Z_{mk}\to \Z_{mkl}$ are embeddings, then $K(\widehat{n_1\cdot n_2})=K(\hat n_2)\cmp K(\hat n_1)$;
	 \item if $\hat n:\Z_m\to \Z_{mk}$ is an embedding, then $\eta_{mk}(\hat n(x))=K(\hat n)(\eta_m(x))$, for each $x\in \Z_m$.
	\end{enumerate}
	
	We first prove 1.
	So we have to prove that for $x,y\in U$, $K(\hat n)(x+y)=K(\hat n)(x)+K(\hat n)(y)$, and that $K(\hat n)(x)=0$ only if $x=0$.
	Take any $x,y\in U$.
	Then
	\[
	 \begin{array}{rl}
	 K(\hat n)(x+y)=&\seq{[n]_p(x_p+y_p):p\in \mathbb P}\\[1mm] =&\seq{[n]_px_p+[n]_py_p:p\in \mathbb P}\\[1mm] =&\seq{[n]_px_p:p\in \mathbb P}+\seq{[n]_py_p:p\in \mathbb P}\\[1mm] =&K(\hat n)(x)+K(\hat n)(y).
	 \end{array}
	\]
	Now suppose that $K(\hat n)(x)=0$ and that $x\neq 0$.
	Since $x\neq 0$, there is a prime $p$ such that $x_p\notin \Z$.
	Hence $x_p=\frac{a}{p^{\al}}$, where $a\in \Z$, $\al$ is a positive integer, and moreover $\gcd(a,p)=1$.
	Since $K(\hat n)(x)=0$, it must be that $[n]_px_p=\frac{[n]_pa}{p^{\al}}\in \Z$.
	But this is not possible because $\gcd([n]_pa,p)=1$.
	Hence $K(\hat n)$ is an injective homomorphism and 1 is proved.
	
	Next, we prove 2.
	Take arbitrary $x\in U$.
	Then
	\[
	 \begin{array}{rl}
	  K(\widehat{n_1\cdot n_2})(x)=&\seq{[n_1n_2]_px_p:p\in \mathbb P}\\[1mm]
	  =&\seq{[n_1]_p[n_2]_px_p:p\in \mathbb P}\\[1mm]
	  =&K(\hat n_2)(\seq{[n_1]_px_p:p\in \mathbb P})\\[1mm]
	  =&K(\hat n_2)(K(\hat n_1)(x)),
	 \end{array}
	\]
	so 2 is proved as well.
    
	Finally we check condition 3.
	Take any $x\in \Z_m$.
	Since $n=k\cdot t$ for $t$ such that $\gcd(t,mk)=1$, it must be that $[n]_p=[kt]_p=t\cdot[k]_p$ whenever $p\mid mk$.
	On the other hand, if $p\nmid mk$, it must be that $\frac{[m]_p}{m}=1$.
	Hence in this case both $\frac{t\cdot [k]_p\cdot x\cdot [m]_p}{m}$ and $\frac{[n]_p\cdot x\cdot [m]_p}{m}$ belong to $\Z$, so $\frac{t\cdot [k]_p\cdot x\cdot [m]_p}{m}\equiv_{U(p)}\frac{[n]_p\cdot x\cdot [m]_p}{m}$. So we have
	\[\displaystyle
	\begin{array}{rl}
	  \eta_{mk}(\hat n(x))=&\eta_{mk}(nx)\\[1mm]
	  =&\seq{\frac{nx\cdot [mk]_p}{mk}:p\in \mathbb P}\\[1mm]
	  =&\seq{\frac{k\cdot t\cdot x\cdot [mk]_p}{mk}:p\in \mathbb P}\\[1mm]
	  =&\seq{\frac{t\cdot [k]_p\cdot x\cdot [m]_p}{m}:p\in \mathbb P}\\[1mm]
	  =&\seq{\frac{[n]_p\cdot x\cdot [m]_p}{m}:p\in \mathbb P}\\[1mm]
	  =&K(\hat n)(\seq{\frac{x\cdot [m]_p}{m}:p\in \mathbb P})\\[1mm]
	  =&K(\hat n)(\eta_m(x)).
	\end{array}
	\]
        This proves 3, and finishes the proof of the theorem.
\end{pf}


\begin{wn}
	Every finite cyclic group is uniformly homogeneous.
\end{wn}

\begin{pf}
	Let $\map K \Cee{\emb{U}}$ be an ultimate \Katetov\ functor, where $U = \Q\by \Z$ is the \fra\ limit of $\Cee$.
	Let $\map{e}{\Z_k}{\Z_n}$ be an embedding.
	We may think of $\Z_n$ as the unique subgroup of $U$ of size $n$. Finally, $E(f) = K(f) \rest \Z_n$ provides an extension operator from $\aut(\Z_k)$ to $\aut(\Z_n)$.
	By Proposition~\ref{PropSDnldvno}, $\Z_n$ is uniformly homogeneous.
\end{pf}




\begin{thebibliography}{99}


\bibitem{Fraisse} {\sc R. \fra}, {\it Sur l'extension aux relations de quelques propri\'et\'es des ordres}, Ann. Sci. Ecole Norm. Sup. (3) {\bf71} (1954) 363--388.

\bibitem{Hodges} {\sc W. Hodges}, {\it  Model Theory}, Encyclopedia of Mathematics and its Applications 42, Cambridge University Press 1993.



\bibitem{KubMas} {\sc W. Kubi\'s, D. Ma\v sulovi\'c}, {\it \Katetov\ functors}, Applied Categorical Structures 25 (2017) 569--602.

\bibitem{Lachlan} {\sc A. H. Lachlan}, {\it Finite homogeneous simple digraphs}. Proceedings of the {H}erbrand symposium ({M}arseilles, 1981), 189--208, Stud. Logic Found. Math., 107, North-Holland, Amsterdam, 1982.

\bibitem{Lang} {\sc S. Lang}, {\it Algebra. Revised third edition}, Graduate Texts in Mathematics, 211. Springer-Verlag, New York, 2002.


\end{thebibliography}
\end{document}